\documentclass[11pt]{article}
\usepackage[english]{babel}
\usepackage{amssymb,amsmath,amsthm}
\newtheorem{theorem}{Theorem}[section]
\newtheorem{lemma}[theorem]{Lemma}
\newtheorem{proposition}[theorem]{Proposition}
\newtheorem{corollary}[theorem]{Corollary}

{\theoremstyle{definition}}
{\theoremstyle{definition}\newtheorem{example}[theorem]{Example}}
{\theoremstyle{definition}}
{\theoremstyle{definition}}

\newtheorem*{thmgs}{Theorem GS}
\newtheorem*{thmgs1}{Theorem GS$'$}
\newtheorem*{cona}{Conjecture A}
\newtheorem*{thma}{Theorem A}
\newtheorem*{conv}{Conjecture V}

\numberwithin{equation}{section}

\def\C{{\mathbb C}}
\def\N{{\mathbb N}}
\def\Z{{\mathbb Z}}
\def\R{{\mathbb R}}
\def\Q{{\mathbb Q}}
\def\K{{\mathbb K}}
\def\F{{\mathcal F}}

\def\E{{\mathcal E}}

\def\epsilon{\varepsilon}
\def\kappa{\varkappa}
\def\phi{\varphi}
\def\leq{\leqslant}
\def\geq{\geqslant}
\def\slim{\mathop{\hbox{$\overline{\hbox{\rm lim}}$}}\limits}
\def\ilim{\mathop{\hbox{$\underline{\hbox{\rm lim}}$}}\limits}

\def\dim{{\rm dim}\,}
\def\ssub#1#2{#1_{{}_{{\scriptstyle #2}}}}
\def\dimk{{\ssub{\dim}{\K}\,}}

\def\spann{\hbox{\tt span}\,}

\def\deg{\hbox{\tt deg}\,}

\def\ootimes{\,{\text{$\scriptstyle\otimes$}}\,}

\title{The Golod-Shafarevich inequality for Hilbert series of quadratic algebras and the Anick
conjecture}

\author{Natalia Iyudu and Stanislav Shkarin}

\date{}

\begin{document}

\maketitle

\begin{abstract}
We study the question on whether the famous Golod-Shafarevich
estimate, which gives a lower bound for the Hilbert series of a
(noncommutative) algebra, is attained. This question was considered
by Anick in his 1983 paper 'Generic algebras and CW-complexes',
Princeton Univ. Press., where he proved that the estimate is
attained for the number of quadratic relations $d \leq \frac{n^2}{4}$ and $d \geq \frac{n^2}{2}$, and conjectured  that it is the case
for any number of quadratic relations.
The particular point where the number of  relations is equal to $ \frac{n(n-1)}{2}$ was addressed by Vershik. He conjectured that a generic algebra with this number of relations is finite dimensional.

We prove that over any infinite field, the Anick conjecture holds for
  $d \geq \frac{4(n^2+n)}{9}$ and arbitrary number of generators, and confirm the Vershik conjecture over any field of characteristic $0$. We give also a series of related asymptotic results.

\end{abstract}




\vspace{5mm}

\section{Introduction}

Let $\F(n,\K)=\K\langle
x_1,\dots,x_n\rangle$ be a free associative algebra
with $n$ generators $x_1,\dots,x_n$, over a field $\K$. Recall that the free algebra
carries the natural degree grading
$$
\F(n,\K)=\bigoplus _{k=0}^\infty \F_k(n,\K),\ \ \text{where}\ \
\F_k(n,\K)=\spann_\K\{x_{j_1}{\dots}x_{j_k}:1\leq j_1,\dots,j_k\leq
n\}.
$$

To define this grading we suppose that generators $x_i$ all have degree one.
We deal with {\it quadratic
algebras generated in degree one}, that is, algebras given by homogeneous relations of
degree 2:
\begin{equation}\label{R}
R=\K\langle x_1,\dots,x_n\rangle\,/\,I,\ \ \text{where}\ \ I={\tt
Id}\,\{f_1,\dots,f_d\}
\end{equation}
is the ideal generated by
\begin{equation}\label{xi}
f_1,\dots,f_d\in \F_2(n,\K):\ \ f_j=\sum_{k,m=1}^n
c_{j,k,m}x_kx_m,\quad c_{j,k,m}\in\K.
\end{equation}
Since relations (\ref{xi}) are homogeneous, an algebra $R$ inherits grading from the free algebra
$\F(n,\K)$:
$$
R=\bigoplus_{k=0}^\infty R_k,\ \ \text{where}\ \ I=\bigoplus_{k=0}^\infty I_k,\ \
\text{for}\ \ I_k=I\cap
\F_k(n,\K)\ \ \text{and}\ \ R_k=\F_k(n,\K)/I_k.
$$
Recall also that the {\it Hilbert series} of $R$ is a polynomial generating function associated to the sequence
of dimensions of graded components
$a_q=\dimk R_q$:

\begin{equation}\label{HI}
 H_R(t)=\sum_{q=0}^\infty (\dimk R_q)\,t^q.
\end{equation}
It belongs to the ring ${\Bbb Z} [[t]]$  of formal power series on one variable and we consider the following ordering on it.
For two power series $a(t)=\sum\limits_{j=0}^\infty a_jt^j$ and
$b(t)=\sum\limits_{j=0}^\infty b_jt^j$ (with real coefficients) we
write $a(t)\geq b(t)$ if $a_j\geq b_j$ for any $j\in\Z_+$. For such a
power series $a(t)=\sum\limits_{j=0}^\infty a_jt^j$, we denote by
$|a(t)|$ a series  obtained from
$a(t)$ by replacing all coefficients starting from the
first negative one by zeros.

The famous result due to
Golod and Shafarevich \cite{gosh} gives a lower bound for the Hilbert series of algebras with $n$ generators and $d$ quadratic relations. (Throughout the text, whenever we are talking on the number of relations, we mean the number of linearly independent relations).

\begin{thmgs}Let $\K$ be a field, $n\in\N$, $0\leq d\leq
n^2$ and $R$ be a quadratic $\K$-algebra with $n$ generators and $d$
relations. Then $H_R(t)\geq |(1-nt+dt^2)^{-1}|$.
\end{thmgs}

Let us note, that  although we formulated above the Golod-Shafarevich estimate only for quadratic algebras, it is known for algebra with any, finite or infinite number of relations. Namely,
it is as follows:

$$ |(1-nt+\sum\limits_{i=2}^{\infty} r_i t^i)^{-1}| \leq H_A,$$

where $r_i$ stands for the number of relations of degree $i$.

This estimate allowed to construct a counterexample for the Kurosh problem on the nilpotency of nil algebra and to the General Burnside Problem on the existence of a finitely generated  infinite torsion $p$-group. It is also recognized due to other applications to $p$-groups  and class field theory \cite{Golod, gosh}.

It will be convenient for our purposes to state the Golod-Safaevich theorem also in terms of the numbers
\begin{equation}\label{hm}
h_q(\K,n,d)=\min_{R\in {\cal R}_{n,d}} \dim R_q,
\end{equation}
where ${\cal R}_{n,d}$ is the set of all quadratic $\K$-algebras $R$
with $n$ generators and $d$ relations. For $n\in\N$ and $0\leq d\leq
n^2$, consider the series
\begin{equation}\label{Hm}
H^{\rm min}_{\K,n,d}(t)=\sum_{q=0}^\infty h_q(\K,n,d)t^q.
\end{equation}
Then Theorem~GS admits the following equivalent form.

\begin{thmgs1}Let $\K$ be a field, $n\in\N$ and $0\leq d\leq
n^2$. Then $H^{\rm min}_{\K,n,d}(t)\geq |(1-nt+dt^2)^{-1}|$.
\end{thmgs1}

Note that {\it a priori} it is not clear why the algebra with the
series $H^{\min}_{\K,n,d}$ should exist in the class $R_{n,d}$. In
fact, it is not difficult to show that not only it does exist, but
it is in 'general position' in one or another sense. Usually by
'generic quadratic algebra' we mean generic in the sense of Zariski
topology. Namely, we consider an algebra from $R_{n,d}$ as a point in
$\K^{n^2d}$, defined by the tuple of all coefficients of its defining
relations. Then we say that a
generic quadratic $\K$-algebra with $n$ relations and $d$ generators
has a property $\cal P$, if the set $\{c_{j,k,m}\}\in \K^{n^2d}$ of
coefficient vectors for which the corresponding algebra $R$ defined
in (\ref{R}) has property $\cal P$ contains a dense Zariski open subset
of $\K^{n^2d}$. The following proposition is a well-known fact, see,
for instance \cite{ani2,cana,popo}.

\begin{proposition}\label{geN} Let $\K$ be an infinite field and
$n\in\N$, $0\leq d\leq n^2$. Then $\dimk R_q=h_q(\K,n,d)$ for a
generic quadratic $\K$-algebra $R$ with $n$ generators and $d$
relations. In particular, if $H^{\rm min}_{\K,n,d}(t)$ is a
polynomial, then $H_R(t)=H^{\rm min}_{\K,n,d}(t)$ for a generic quadratic
$\K$-algebra $R$ with $n$ generators and $d$ relations.
\end{proposition}

In the case  $H^{\min}_{\K,n,d}$ is not a polynomial, there are
more subtleties involved in the question whether a generic algebra
has this series. There are arguments (see \cite{ufnar}) showing that this
is the case, when $\K$ is an uncountable algebraically closed field.
In \cite{cana} we suggested a slightly modified notion of a
'generic' algebra, in the case $\K=\R$.
Namely, we say that the generic in the Lebesgue sense algebra from
$R_{n,d}$ has the property $\cal P$ if the set of algebras not having $\cal P$
has Lebesgue measure zero. We show that in this (weaker) sense a
generic algebra has the series $H^{\min}$  even if it is infinite.

In his 1983 paper "Generic algebras and CW-complexes",
Princeton Univ. Press, \cite{ani1}, Anick studied the behavior of Hilbert
series of algebras given by relations and
formulated the following conjecture.

\begin{cona}For any infinite field $\K$, any $n,q\in\N$ and
$0\leq d\leq n^2$, a generic quadratic $\K$-algebra $R$ with $n$
generators and $d$ relations $\dim R_q$ equals to the $q^{\rm th}$
coefficient of the series $|(1-nt+dt^2)^{-1}|$. Equivalently,
$H^{\rm min}_{\K,n,d}(t)=|(1-nt+dt^2)^{-1}|$.
\end{cona}

In other words, Conjecture~A states that the lower estimate of the
Hilbert series by Golod and Shafarevich is attained  and a generic algebra has the minimal Hilbert series.

This question is very important to clarify in the light of key problems in the ring theory concerned with the behavior of nilpotent elements. Examples of such problems are the problem on the existence of simple nil ring,
solved in affirmative by A.Smoktuniwicz \cite{Ag8}, the K\"othe conjecture, the Burnside type problem for finitely presented rings (see \cite{Ag9}).

Values  of terms $h_q(\K,n,d)$ for $q=0,1,2$ are obvious for an
arbitrary algebra given by $n$ generators and $d$ relations: $h_0(\K,n,d)=1$, $h_1(\K,n,d)=n$ and
$h_2(\K,n,d)=n^2-d$. Anick proved \cite{ani1,popo} that his conjecture holds also for
$q=3$.


\begin{thma}Let $\K$ be any field, $n\in\N$ and $0\leq d\leq n^2$.
Then
\begin{equation}\label{com3}
h_3(\K,n,d)=\left\{\begin{array}{ll}0&\text{if\ \ $d\geq \frac{n^2}2$};\\
n^3-2nd&\text{if\ \ $d<\frac{n^2}2$.}\end{array}\right.
\end{equation}
\end{thma}

Since the number in the right-hand side of (\ref{com3}) happens to
coincide with the third coefficient near $t^3$ in $|(1-nt+dt^2)^{-1}|$,
 Theorem~A proves Conjecture~A in the case $q=3$ and in the
case $d\geq \frac{n^2}2$. Conjecture~A is also known to be true if
$d\leq \frac{n^2}{4}$ \cite{ufnar,popo}. The region
$\frac{n^2}4<d<\frac{n^2}2$ remained a white zone so far.
 Let us note that for $d>\frac{n^2}4$, the series
$|(1-nt+dt^2)^{-1}|$ is a polynomial. Thus Conjecture~A, if true,
implies that a generic quadratic $\K$-algebra with infinite $\K$,
$n$ generators and $d>\frac{n^2}4$ relations is finite dimensional.

In  \cite{vers} Vershik formulated a conjecture, which addresses a specific point
of the 'difficult interval' $\frac{n^2}4<d<\frac{n^2}2$, $d=\frac{n(n-1)}2$, which is
 the number of relations defining the algebra of
commutative polynomials or any PBW algebra.

\begin{conv} Let $n\in\N$, $n\geq 3$. Then a generic quadratic
$\C$-algebra with $n$ generators and $\frac{n(n-1)}{2}$ relations is
finite dimensional.
\end{conv}

As it is mentioned in \cite{popo} there was an attempt to prove this conjecture in
\cite{zvy}, but the argument there was incorrect.

Our goal in this paper is to move the frame of the interval
$(\frac{n^2}{4},\frac{n^2}{2})$ which remains unknown since  the
Anick's 1983 paper. Namely, we prove the following.

\vspace{5mm}

\begin{theorem}\label{main0} For  any infinite field, the Golod--Shafarevich estimate is
attained for a generic quadratic algebra with $n$ generators and
$d\geq \frac{4(n^2+n)}{9}$ quadratic relations.

Namely, the Hilbert series of the generic algebra is:

$$H(t)=|(1-nt+dt^2)^{-1}|=1+nt+(n^2-d)t^2+(n^3-2nd)t^3.$$\rm

\end{theorem}

The point $d=\frac{n(n-1)}{2}$ falls into the interval from the Theorem\ref{main0}, for big
enough $n$, so we automatically get as a consequence an affirmative
answer to the Vershik's question for $n\geq 17$. After some
additional considerations, we get an affirmative answer for the
Vershik's question for each $n\geq 3$ over a field of
characteristic $0$:

\vspace{5mm}

\begin{theorem}\label{solv}
Let $\K$ be a field of characteristic $0$ and $n\in\N$, $n\geq 3$.
Then a generic quadratic $\K$-algebra $R$ with $n$ relations and
$\frac{n(n-1)}{2}$ relations has the Hilbert series and the
dimension given by the following formula
\begin{equation}\label{hsver}
H_R(t)=\left\{
\begin{array}{ll}
1+nt+\frac{n(n+1)}2 t^2+n^2t^3&\text{if\ \ $n\geq 5;$}\\
1+4t+10t^2+16t^3+t^4&\text{if\ \ $n=4;$}\\
1+3t+6t^2+9t^3+9t^4&\text{if\ \ $n=3$,}
\end{array}\right.
\quad\dimk R=
\left\{\begin{array}{ll}\frac{3n(n+1)+2}2&\text{if\ \ $n\geq 5;$}\\
32&\text{if\ \ $n=4;$}\\ 28&\text{if\ \ $n=3$.}
\end{array}\right.
\end{equation}
\end{theorem}

We will formulate more explicit results later in
the text. In order to illustrate them, we present their asymptotic
versions straight away.

Let us note that the series of related questions on asymptotic characteristics of the Golod-Shafarevich inequality were considered in \cite{Sch}

For a field $\K$, $n,q\in\N$ with $q\geq 2$,
we denote
\begin{equation}\label{dnq}
d(\K,n,q)=\min\{d\in\N:h_q(\K,n,d)=0\}.
\end{equation}
That is, $d(\K,n,q)$ is the minimal $d$ for which there is a
quadratic $\K$-algebra $R$ with $n$ generators and $d$ relations
satisfying $R_q=\{0\}$. Obviously, $d(\K,n,2)=n^2$. Similarly
\begin{equation}\label{dni}
d(\K,n,\infty)=\min\Bigl\{d\in\N: \min_{q\in\N}h_q(\K,n,d)=0\Bigr\}.
\end{equation}
That is, $d(\K,n,\infty)$ is the minimal $d$ for which there is a
finite dimensional quadratic $\K$-algebra $R$ with $n$ generators
and $d$ relations. In order to formulate our asymptotic results we
need the following lemma.

\begin{lemma}\label{limi} Let $\K$ be a field and $q\in\N$, $q\geq
2$ or $q=\infty$. Then the limit
$\lim\limits_{n\to\infty}\frac{d(\K,n,q)}{n^2}=\alpha(\K,q)$ does
exist and
\begin{equation}\label{liinf}
\alpha(\K,q)=\lim\limits_{n\to\infty}\frac{d(\K,n,q)}{n^2}=\inf\biggl\{
\frac{d(\K,n,q)}{n^2}:n\in\N\biggr\}.
\end{equation}
Moreover, $\{\alpha(\K,q)\}_{q\geq 3}$ is decreasing,
$\alpha(\K,\infty)=\lim\limits_{q\to\infty}\alpha(\K,q)\geq\frac14$
and $\alpha(\K,3)=\frac12$.
\end{lemma}

\begin{theorem}\label{lith1}The equalities
$\alpha(\K,4)=\frac{3-\sqrt5}{2}$ and $\alpha(\K,5)=\frac13$ hold
for any infinite field. Moreover, $\frac14\leq
\alpha(\K,\infty)\leq\alpha(\K,6)\leq \frac{5}{16}$ for any field
$\K$ of characteristic $0$.
\end{theorem}

\begin{corollary}\label{lith2} Let $\K$ be an infinite field
and $\ilim\limits_{n\to\infty}\frac{d_n}{n^2}>\frac{3-\sqrt5}{2}$
with $n,d_n\in\N$ and $d_n\leq n^2$. Then for any sufficiently large
$n$, a generic quadratic $\K$-algebra with $n$ generators and $d_n$
relations has Hilbert series
$1+nt+(n^2-d_n)t^2+\max\{0,(n^3-2d_nn^2)\}t^3=|(1-nt+d_nt^2)^{-1}|$.
\end{corollary}

\begin{proof}Indeed, by Theorem~\ref{lith1}, the Hilbert series in question
must be a polynomial of degree at most $3$, whose specific shape is
determined by Theorem~A and Proposition~\ref{geN}.
\end{proof}

\begin{corollary}\label{lith2a} Let $\K$ be an infinite field
and $\ilim\limits_{n\to\infty}\frac{d_n}{n^2}>\frac{1}{3}$ with
$n,d_n\in\N$ and $d_n\leq n^2$. Then for any sufficiently large $n$,
the Hilbert series of a generic quadratic $\K$-algebra with $n$
generators and $d_n$ relations is a polynomial of degree at most
$4$.
\end{corollary}

\begin{corollary}\label{lith3} Let $\K$ be a field of characteristic
$0$ and $\ilim\limits_{n\to\infty}\frac{d_n}{n^2}>\frac{5}{16}$ with
$n,d_n\in\N$ and $d_n\leq n^2$. Then for any sufficiently large $n$,
a generic quadratic $\K$-algebra with $n$ generators and $d_n$
relations has Hilbert series being a polynomial of degree at most
$5$ and therefore is finite dimensional.
\end{corollary}

\section{Notations and preliminary facts \label{genob}}

Let $\K$ be a field, $n,d,q\in\N$, $1\leq d\leq n^2$, $q\geq 3$ and
$\{c_{j,k,m}:1\leq j\leq d,1\leq k,m\leq n\}$ be variables taking
values in $\K$. Consider the ideal $I_c$ with $c=\{c_{j,k,m}\}$ in
$\K\langle x_1,\dots,x_n\rangle$ generated by $f_1,\dots,f_d$,
$f_j=\sum\limits_{k,m=1}^n c_{j,k,m}x_kx_m$ and the algebra
$R_c=\K\langle x_1,\dots,x_n\rangle/I_c$. Clearly the $q^{\rm th}$
homogeneous component $(I_c)_q$ is spanned by $\mu f_j\nu$, where
$1\leq j\leq d$ and $\mu,\nu$ are two monomials in $\K\langle
x_1,\dots,x_n\rangle$ with $\deg \mu+\deg \nu=q-2$. Hence $(I_c)_q$
is the image of the linear operator $L_c:\K^{\Omega}\to \F_q(n,\K)$,
where $\Omega$ is the set of triples $(j,\mu,\nu)$ with $1\leq j\leq
d$, $\mu$ and $\nu$ being monomials satisfying $\deg \mu+\deg
\nu=q-2$ and $L_c$ sends the standard basic vector $e_{j,\mu,\nu}$
to $\mu f_j\nu$. Then the dimension of $(I_c)_q$ equals to the rank
${\tt rk}\,L_c$ of $L_c$. Hence $\dim (R_c)_q=n^q-\dim
(I_c)_q=n^q-r_c$, where $r_c={\tt rk}\,L_c$. It immediately follows
that
$$
h_q(\K,n,d)=n^q-r,\ \ \text{where}\ \ r=\max_c r_c.
$$
Since the rank of a linear map equals to the maximal size of its
square submatrix with non-zero determinant, there exist an $r\times
r$ submatrix of the rectangular matrix of $L_c$, whose determinant
$\delta(c)$ is non-zero for some $c$. On the other hand, $\delta(c)$
is a polynomial in $c_{j,k,m}$ with integer coefficients. Since a
polynomial with integer coefficients over a field of characteristic
$0$ defines a zero function if and only if all its coefficients are
zero, we see that the numbers $h_q(\K,n,d)$ do not depend on the
choice of $\K$ provided $\K$ has characteristic $0$.

Now if $p$ is a prime number and $\K=\Z_p$, then the fact that
$\delta(c)$ is non-zero for some $c$ implies that the coefficients
of $\delta(c)$ are not all zeros as elements of $\Z_p$. Hence some
of the coefficients of $\delta(c)$ considered as a polynomial with
coefficients in $\Z$ are not multiples of $p$ and therefore are
non-zero. Hence $\delta(c)$ remains non-zero for some $c$ after
replacing the field $\Z_p$ by $\Q$. It follows that $h_q(\Q,n,d)\leq
h_q(\Z_p,n,d)$. Similar argument shows that if $h_q(\K,n,d)$ does
not depend on the choice of an infinite field $\K$ of a fixed
positive characteristic $p$ and that $h_q(\K,n,d)\leq h_q(\Z_p,n,d)$
for any field $\K$ of characteristic $p$.

Next, if $R$ is a quadratic $\K$-algebra with $n$ generators
$x_1,\dots,x_n$ and $d$ relations, then the quotient $R_0$ of $R$ by
the ideal generated by $x_{n'+1},\dots,x_n$ is a quadratic
$\K$-algebra with $n'$ generators and $d$ relations, whose
homogeneous components are quotients of the homogeneous components
of $R$. It follows that $h_q(\K,n',d)\leq h_q(\K,n,d)$ if $n'\leq
n$. On the other hand adding new relations to an algebra can only
decrease the dimension of its components. Hence $h_q(\K,n,d)\geq
h_q(\K,n,d')$ if $d\leq d'$. Above observations are summarized in
the following proposition.

\begin{proposition}\label{fiel} For any field $\K$,
$h_q(\K,n,d)$ increase with respect to $n$ and decrease with respect
to $d$. Moreover, if $n\in\N$, $q,d\in\Z_+$ and $0\leq d\leq n^2$
and $p$ is a prime number, then $H^{\rm min}_{\K,n,d}(t)=H^{\rm
min}_{\Q,n,d}(t)\leq H^{\rm min}_{\Z_p,n,d}(t)$ for any field $\K$
of characteristic zero and $H^{\rm min}_{\K_1,n,d}(t)=H^{\rm
min}_{\K_2,n,d}(t)\leq H^{\rm min}_{\Z_p,n,d}(t)$ for any two
infinite fields $\K_1$ and $\K_2$ of characteristic $p$.
\end{proposition}

In what follows it is convenient to give an alternative definition
of the numbers $h_q(\K,n,d)$. We need the following notation. Let
$E$ be a vector space over a field $\K$. For $k\in\Z_+$, we denote
the $k^{\rm th}$ tensor power of $E$ by $E^{\otimes k}$. That is,
$E^{\otimes 0}=\K$, $E^{\otimes 1}=E$ and $E^{\otimes k}=E\otimes
{\dots}\otimes E$ is the tensor product of $k$ copies of $E$. If $L$
is a subspace of $E^{\otimes 2}=E\otimes E$, then for $k\geq 2$ we
can define the subspaces $\E_k(L,E)$ of $E^{\otimes k}$ inductively:
$\E_2(L,E)=L$ and $\E_{k+1}(L,E)=E\otimes \E_k(L,E)\cap
\E_k(L,E)\otimes E$. We can also write two explicit expressions for
the space $\E_k(L,E)$:
\begin{align}\label{ekl1}
\E_k(L,E)&=\left\{\begin{array}{ll}(E\otimes
L^{\otimes(k-2)/2}\otimes E)\cap L^{\otimes(k/2)}&\text{if $k$ is
even;}\\(E\otimes L^{\otimes(k-1)/2})\cap (L^{\otimes(k-1)/2}\otimes
E)&\text{if $k$ is odd,}
\end{array}
\right.
\\
\label{ekl} \E_k(L,E)&=\bigcap_{j=1}^{k-1}L^{k,j},\quad
\text{where}\ \ L^{k,j}=E^{\otimes (j-1)}\otimes L\otimes E^{\otimes
(k-1-j)}\ \ \text{for $1\leq j\leq k-1$}.
\end{align}

If $E$ is an $n$-dimensional vector space over $\K$ with a fixed
basis $\{e_1,\dots,e_n\}$, we consider the symmetric bilinear form
$[\cdot,\cdot]_j:E^{\otimes j}\times E^{\otimes j}\to\K$ such that
$$
[e_{m_1}\ootimes{\dots}\ootimes
e_{m_j},e_{r_1}\ootimes{\dots}\ootimes e_{r_j}]_j=\delta_{m,r},\ \
\text{where $m,r\in\{1,\dots,n\}^j$.}
$$
For a subspace $N$ of $E^{\otimes j}$ we write
$$
N^\perp=\{f\in E^{\otimes j}:[\eta,f]_j=0\ \ \text{for all}\ \
\eta\in N\}.
$$
From the definition of $[\cdot,\cdot]_j$ it easily follows that the
space $N^\perp$ is naturally isomorphic to the space of linear
functionals on $E^{\otimes j}$ annihilating $N$. Hence
\begin{equation}\label{dimp}
\dim N+\dim N^\perp=\dim E^{\otimes j}=n^j.
\end{equation}
Moreover, one can easily see that $(N^\perp)^\perp=N$ and therefore
the set of $N^\perp$ for all $d$-dimensional subspaces $N$ of
$E^{\otimes j}$ coincides with the set of all $(n^j-d)$-dimensional
subspaces of $E^{\otimes j}$.

\begin{lemma}\label{rel} Let $\K$ be a field, $n\in\N$,
$1\leq d\leq n^2$, $E$ be an $n$-dimensional vector space over $\K$
with a basis $\{e_1,\dots,e_n\}$ and $R$ be a quadratic $\K$-algebra
with $n$ generators $x_1,\dots,x_n$ and $d$
relations $f_s=\!\!\!\!\sum\limits_{1\leq a,b\leq n}\!\!\!
c_{s,a,b}x_ax_b$. Let also $M$ be the $d$-dimensional subspace of
$E\otimes E$ spanned by $\sum\limits_{1\leq a,b\leq n}\!\!\!
c_{s,a,b}e_a\ootimes e_b$ for $1\leq s\leq d$. Then $\dim
R_q=\dim\E_q(M^\perp,E)$ for any $q\geq 2$.
\end{lemma}

\begin{proof} Let $q\geq2$. Under the linear isomorphism between
${\cal F}_q(\K,n)$ which sends $x_{m_1}{\dots}x_{m_q}$ to
$e_{m_1}\ootimes{\dots}\ootimes e_{m_q}$, the $q^{\rm th}$
homogeneous component $I_q$ of the ideal $I$ generated by
$f_1,\dots,f_d$ is mapped onto the subspace
$$
{\cal M}=\sum_{j=1}^{q-1} M^{q,j}
$$
of $E^{\otimes q}$, where $M^{q,j}$ are defined in (\ref{ekl}).
Using the last display, is straightforward to see that
$$
{\cal
M}^\perp=\bigcap_{j=1}^{q-1}(M^{q,j})^\perp=\bigcap_{j=1}^{q-1}(M^\perp)^{q,j}=
\E_q(M^\perp,E),
$$
where the latter space is defined in (\ref{ekl}). Now using
(\ref{dimp}), we see that
$$
\dim R_q=n^q-\dim I_q= n^q-\dim{\cal M}=\dim\E_q(M^\perp,E),
$$
which completes the proof. \end{proof}

Next lemma relates the spaces $\E_q(L,E)$ and the numbers
$h_q(\K,n,d)$.

\begin{lemma}\label{tepr} Let $\K$ be a field, $q,n\in\N$, $d\in\Z_+$,
$q\geq 3$, $0\leq d\leq n^2$. Then
\begin{equation}\label{tp1}
h_q(\K,n,d)=\min\{\dim\E_q(L,E):\dim L=n^2-d\},
\end{equation}
where the minimum is taken over $(n^2-d)$-dimensional subspaces $L$
of $E\otimes E$ with $E$ being an $n$-dimensional vector space over
$\K$.
\end{lemma}

\begin{proof} Follows immediately from the above lemma and the fact
that the map $M\mapsto M^\perp$ is a bijection between the sets of
$d$-dimensional and $(n^q-d)$-dimensional subspaces of $E^{\otimes
q}$.
\end{proof}

\section{Main lemma}

\begin{lemma}\label{main}
Let $\K$ be a field, $n,q,m\in\N$, $q\geq 2$, $E$ be an
$n$-dimensional vector space and
$$
E_1\subset E_2\subset{\dots}\subset E_{n-1}\subset E_n=E
$$
is an increasing chain of subsets of $E$ such that $\dim E_j=j$ for
$1\leq j\leq n$. Let also $L$ be a subspace of $E\otimes E$,
$n_1,\dots,n_m\in\{1,\dots,n\}$ and
$$
G=\bigoplus_{j=1}^m E_{n_j}\ \ \ \text{and}\ \ \ L_G= \bigoplus
_{1\leq j,k\leq m} L_{j,k},\ \ \text{where}\ \ L_{j,k}=(E_j\otimes
E_k)\cap L.
$$
Then
\begin{equation}\label{DIM}
\dim G=\sum_{j=1}^m n_j,\ \ \dim L_G=\sum_{j,k=1}^m \dim
L_{n_j,n_k}\ \ \text{and}\ \ \dim \E_q(L_G,G)\leq m^q\dim \E_q(L,E).
\end{equation}
\end{lemma}

\begin{proof} The two equalities in (\ref{DIM}) are trivial.
Using the natural decomposition
\begin{align*}
&G^{\otimes q}=\bigoplus_{a_1,\dots,a_q=1}^m
E_{n_{a_1}}\otimes{\dots}\otimes E_{n_{a_q}},\ \ \text{we see that}\
\ \E_q(L_G)=\bigoplus_{a_1,\dots,a_q=1}^m F_a,\ \text{where}
\\
F_a&=(L_{n_{a_1},n_{a_2}}\otimes E_{n_{a_3}}\otimes{\dots}\otimes
E_{n_{a_q}})\cap (E_{n_{a_1}}\otimes L_{n_{a_2},n_{a_3}}\otimes
E_{n_{a_4}}\otimes {\dots}\otimes E_{n_{a_q}})\cap {\dots}
\\
&\qquad\qquad \dots\cap (E_{n_{a_1}}\otimes {\dots}\otimes
E_{n_{a_{q-2}}}\otimes L_{n_{a_{q-1}},n_{a_{q}}}).
\end{align*}
Clearly each $F_a$ is isomorphic to a subspace of $\E_q(L,E)$ and
therefore $\dim F_a\leq \dim \E_q(L,E)$ for any $a$. Now since
$\E_q(L_G,G)$ is the sum of $F_a$ and there are $m^q$ multi-indices
$a$, we get $\dim \E_q(L_G,G)\leq m^q\dim \E_q(L,E)$.
\end{proof}

\begin{lemma}\label{main1}
Let $\K$ be a field, $n,m\in\N$, $q\geq 2$ and $d\in\Z_+$, $0\leq
d\leq n^2$. Then
\begin{equation}\label{ineq1}
h_{q}(\K,mn,m^2d)\leq m^qh_q(\K,n,d).
\end{equation}
\end{lemma}

\begin{proof} Let $h=h_q(\K,n,d)$ and $E$ be an $n$-dimensional
vector space over $\K$. By Lemma~\ref{tepr}, there exists an
$(n^2-d)$-dimensional subspace $L$ of $E\otimes E$ such that $\dim
\E_q(L)=h$. Let $G$ be the direct sum of $m$ copies of $E$. Clearly
$\dim G=nm$. Applying Lemma~\ref{main} with $n_1={\dots}=n_m=n$, we
find a subspace $L_G$ of $G\otimes G$ of dimension
$m^2(n^2-d)=(mn)^2-m^2d$ such that $\dim \E(L_G,G)\leq m^qh$. By
Lemma~\ref{tepr}, $h_{q}(\K,mn,m^2d)\leq m^qh$, which is the desired
inequality.
\end{proof}

\begin{corollary}\label{main2}
Let $\K$ be a field, $n,m\in\N$, $q\geq 3$ and $d\in\Z_+$, $0\leq
d\leq n^2$. If $h_q(\K,n,d)=0$, then $h_{q}(\K,mn,m^2d)=0$.
\end{corollary}

\begin{corollary}\label{main3}
Let $\K$ be a field, $n,m\in\N$ and $d\in\Z_+$, $0\leq d\leq n^2$.
If $H^{\rm min}_{\K,n,d}(t)=|(1-nt+dt^2)^{-1}|$, then $H^{\rm
min}_{\K,nm,m^2d}(t)=|(1-nmt+dm^2t^2)^{-1}|$.
\end{corollary}

\begin{proof} Lemma~\ref{main1} implies that $H^{\rm
min}_{\K,nm,m^2d}(t)\leq H^{\rm
min}_{\K,n,d}(mt)=|(1-nmt+dm^2t^2)^{-1}|$. The opposite inequality
follows from Theorem~GS.
\end{proof}

\subsection{Proof of Lemma~\ref{limi}}

Let $\K$ be a field and $3\leq q\leq \infty$. For each $k\in\N$, let
$d_k=d(\K,k,q)$ be the numbers defined by (\ref{dnq}) and
(\ref{dni}). By definition of $d_n$, $h_r(\K,n,d_n)=0$ for some
$r\in\N$, $3\leq r\leq q$ (actually $r=q$ if $q<\infty$). Fix
$n\in\N$ and let $k\in\N$. Then there exists $m\in\N$ and
$j\in\{0,\dots,n-1\}$ such that $k=mn-j$. By Proposition~\ref{geN}
and Corollary~\ref{main2}, $h_r(k,m^2d_n)\leq h_r(\K,mn,m^2d_n)=0$.
Hence $h_r(k,m^2d_n)=0$ and therefore $d_k\leq m^2d_n$. Since
$k=mn-j$ and $j\leq n-1$, we have $m\leq \frac{k+n-1}{n}$. Thus
$d_k\leq d_n\frac{(k+n-1)^2}{n^2}$ for any $k\in\N$. Equivalently,
$$
\frac{d_k}{k^2}\leq \frac{d_n(k+n-1)^2}{n^2k^2}\ \ \text{for any}\ \
k,n\in\N.
$$
Passing to the limit as $k\to\infty$, we get
$\slim\limits_{k\to\infty}\frac{d_k}{k^2}\leq \frac{d_n}{n^2}$ for
any $n\in\N$. Hence
$$
\inf_{n\in\N}\frac{d_n}{n^2}\leq
\ilim_{n\to\infty}\frac{d_n}{n^2}\leq
\slim_{n\to\infty}\frac{d_n}{n^2}\leq \inf_{n\in\N}\frac{d_n}{n^2}.
$$
That is, the limit $\lim\limits_{n\to\infty}\frac{d_n}{n^2}$ does
exist and equals $\inf\limits_{n\in\N}\frac{d_n}{n^2}$.

From the definition of $d(\K,n,q)$ and Proposition~\ref{geN} it
immediately follows that $d(\K,n,q_1)\leq d(\K,n,q_2)$ if $q_1\geq
q_2$. Hence the sequence $\{\alpha(\K,q)\}_{q\geq 3}$ is decreasing.
The equality
$\alpha(\K,\infty)=\lim\limits_{q\to\infty}\alpha(\K,q)$ is also
straightforward. The inequalities $\alpha(\K,3)\geq \frac12$ and
$\alpha(\K,\infty)\geq \frac14$ follow from Theorem~GS. By
Theorem~A, $\alpha(\K,3)\leq\frac12$ and therefore
$\alpha(\K,3)=\frac12$.

\section{Hilbert series of degrees 3 and 4}

In this section we compute $\alpha(\K,4)$ and $\alpha(\K,5)$ for any
infinite field $\K$. We prove certain non-asymptotic estimates and
use them to prove Theorem~\ref{solv}.

\subsection{$\alpha(\K,4)$ for any infinite field}

We use the same idea as in the proof of Lemma~\ref{main}. First, we
need the following lemma.

\begin{lemma}\label{alp4} Let $\K$ be any field, $E$ be a vector
space over $\K$ of dimension $n\in\N$, and $r\in\N$ be such that
$1\leq r<n$ and $d^2+n^2d\leq n^4$, where $d=n^2-r^2$. Then there
exist two subspaces $L$ and $M$ of $E\otimes E$ such that $(L\otimes
L)\cap (E\otimes M\otimes E)=\{0\}$ and $\dim L=\dim M=d$.
\end{lemma}

\begin{proof} Pick a basis $\{e_1,\dots,e_n\}$ in $E$.
Now we consider $M=\spann\{e_j\ootimes e_k:\max\{j,k\}>r\}$ and
$L_0=\spann(\{e_j\ootimes e_k:1\leq j,k\leq r\}\cup\{e_j\ootimes
e_k+e_k\ootimes e_j:1\leq j\leq r<k\leq n\})$. Clearly $M$ and $L_0$
are subspaces of $E\otimes E$ of dimensions $d=n^2-r^2$ and $rn$
respectively. The inequality $d^2+n^2d\leq n^4$ implies that $rn\geq
d$ and therefore $\dim L_0\geq d$. It remains to show that
$(L_0\otimes L_0)\cap (E\otimes M\otimes E)=\{0\}$. Indeed, then any
$d$-dimensional subspace $L$ of $L_0$ satisfies $(L\otimes L)\cap
(E\otimes M\otimes E)=\{0\}$.

Let $\xi\in (L_0\otimes L_0)\cap (E\otimes M\otimes E)$. According
to the definitions of $L_0$ and $M$,
\begin{align*}
\xi&=\sum_{r<\max\{k,l\}\leq n\atop 1\leq j,m\leq n}\!\!\!\!\!
a_{j,k,l,m}e_j\ootimes e_k\ootimes e_l\ootimes e_m=
\\
&=\sum_{1\leq j,k,l,m\leq r}\!\!\!\!\! b_{j,k,l,m}e_j\ootimes
e_k\ootimes e_l\ootimes e_m+\!\!\!\!\!\sum_{1\leq j,k\leq r\atop
1\leq l\leq r<m\leq n}\!\!\!\!\! d_{j,k,l,m}(e_j\ootimes e_k\ootimes
e_l\ootimes e_m+e_j\ootimes e_k\ootimes e_m\ootimes e_l)+
\\
&\quad+\!\!\!\!\!\sum_{1\leq l,m\leq r\atop 1\leq j\leq r<k\leq
n}\!\!\!\!\! s_{j,k,l,m}(e_j\ootimes e_k\ootimes e_l\ootimes
e_m+e_k\ootimes e_j\otimes e_l\ootimes e_m)+
\\
&\quad+\!\!\!\!\!\sum_{1\leq j\leq r<k\leq n\atop 1\leq l\leq
r<m\leq n}\!\!\!\! v_{j,k,l,m}(e_j\ootimes e_k\ootimes e_l\ootimes
e_m+e_j\ootimes e_k\ootimes e_m\ootimes e_l+e_k\ootimes e_j\ootimes
e_l\ootimes e_m+e_k\ootimes e_j\ootimes e_m\ootimes e_l),
\end{align*}
where $a_{j,k,l,m}$, $b_{j,k,l,m}$, $c_{j,k,l,m}$, $d_{j,k,l,m}$,
$s_{j,k,l,m}$ and $v_{j,k,l,m}$ are coefficients from $\K$.

If $j,k,l,m\leq r$, then the basic vector $e_j\ootimes e_k\ootimes
e_l\ootimes e_m$ appears in the above display only once and with the
coefficient $b_{j,k,l,m}$. Hence $b_{j,k,l,m}=0$ for any $j,k,l,m$.
If $j>r$ and $k,l,m\leq r$, then the basic vector $e_j\ootimes
e_k\ootimes e_l\ootimes e_m$ appears in the above display only once
and with the coefficient $s_{k,j,l,m}$. Hence $s_{j,k,l,m}=0$ for
any $j,k,l,m$. If $m>r$ and $j,k,l\leq r$, then the basic vector
$e_j\ootimes e_k\ootimes e_l\ootimes e_m$ appears in the above
display only once and with the coefficient $d_{j,k,l,m}$. Hence
$d_{j,k,l,m}=0$ for any $j,k,l,m$. If $j,m>r$ and $k,l\leq r$, then
the basic vector $e_j\ootimes e_k\ootimes e_l\ootimes e_m$ appears
in the above display only once and with the coefficient
$d_{k,j,l,m}$. Hence $d_{j,k,l,m}=0$ for any $j,k,l,m$. Thus the
right-hand side of the above display vanishes and so does $\xi$.
Hence $(L_0\otimes L_0)\cap (E\otimes M\otimes E)=\{0\}$.
\end{proof}

Just the same way as we speak of generic quadratic algebras with
fixed number of relations and generators, we can speak of generic
vector subspaces of given dimension in a fixed finite dimensional
vector space over an infinite field. Let $\K$ be an infinite field
and $E$ be an $n$-dimensional vector space over $\K$. Using the
argument exactly as in Section~\ref{genob}, one can easily show that
if there exist $d$-dimensional subspaces $L_0$ and $M_0$ of
$E\otimes E$ satisfying $(L_0\otimes L_0)\cap (E\otimes M_0\otimes
E)=\{0\}$, then the equality $(L\otimes L)\cap (E\otimes M\otimes
E)=\{0\}$ holds for generic $d$-dimensional subspaces $L$ and $M$ of
$E\otimes E$. Similarly, if there exist $d$-dimensional subspaces
$L_0$, $N_0$ and $M_0$ of $E\otimes E$ satisfying $(L_0\otimes
N_0)\cap (E\otimes M_0\otimes E)=\{0\}$, then the equality
$(L\otimes N)\cap (E\otimes M\otimes E)=\{0\}$ holds for generic
$d$-dimensional subspaces $L$, $N$ and $M$ of $E\otimes E$. Thus
Lemma~\ref{alp4} implies the following fact.

\begin{lemma}\label{alp4a} Let $\K$ be an infinite field, $E$ be a vector
space over $\K$ of dimension $n\in\N$, and $r\in\N$ be such that
$1\leq r<n$ and $d^2+n^2d\leq n^4$, where $d=n^2-r^2$. Then
$(L\otimes L)\cap (E\otimes M\otimes E)=\{0\}$ and $(L\otimes N)\cap
(E\otimes M\otimes E)=\{0\}$ for generic $d$-dimensional subspaces
$L$, $N$ and $M$ of $E\otimes E$.
\end{lemma}

\begin{proposition}\label{whatnot} There exists a positive constant
$C$ such that
\begin{equation}\label{ine}
\frac{3-\sqrt5}2 n^2<d(\K,n,4)\leq \frac{3-\sqrt5}2 n^2+ Cn^{3/2}
\end{equation}
for any $n\in\N$ and any infinite field $\K$. In particular,
$\alpha(\K,4)=\frac{3-\sqrt5}2$ for any infinite field $\K$.
\end{proposition}

\begin{proof} Let $n,d\in\N$ and $d\leq n^2$. Since the coefficient in
$|(1-nt+dt^2)|$ in front of $t^4$ is $\max\{0,n^4-3n^2d+d^2\}$, it
is positive if $d<\frac{3-\sqrt5}2n^2$. By Theorem~GS,
$h_4(\K,n,d)>0$ if $d<\frac{3-\sqrt5}2n^2$. Hence
$d(\K,n,4)>\frac{3-\sqrt5}2n^2$ for any field $\K$.

Let now $\K$ be an infinite field and $n\in\N$. Choose $m\in\N$ such
that $(m-1)^2<n\leq m^2$. Now let $r$ be the unique integer such
that $\frac{\sqrt5-1}{2}m<r<\frac{\sqrt5-1}{2}m+1$. Since
$r^2>\frac{3-\sqrt5}2m^2$, we have $d<\frac{\sqrt5-1}2m^2$, where
$d=m^2-r^2$. The latter inequality implies $d^2+m^2d<m^4$. Let
$E_1,\dots,E_m$ be $m$-dimensional vector spaces over $\K$ and
$E=E_1\oplus{\dots}\oplus E_m$. Obviously $\dim E=m^2\geq n$.
Consider the space
$$
L=\bigoplus _{j,k=1}^m L_{j,k},
$$
where $L_{j,k}$ is a $d$-dimensional subspace of $E_j\otimes E_k$ if
$j\neq k$ and $L_{j,j}=\{0\}$ for $1\leq j\leq m$. Clearly $\dim
L=(m^2-m)d$. According to (\ref{ekl1}), $\E_4(L,E)=(L\otimes L)\cap
(E\otimes L\otimes E)$. Hence
$$
\E_4(L,E)=\bigoplus_{j,k,l,s=1}^m M_{j,k,l,s},\ \ \text{where}\ \
M_{j,k,l,s}=(L_{j,k}\otimes L_{l,s})\cap (E_j\otimes L_{k,l}\otimes
E_s).
$$
Since $L_{j,j}=0$, we have $M_{j,k,l,s}=\{0\}$ if $j=k$, or $k=l$,
or $l=s$. If $j\neq k$, $k\neq l$ and $l\neq s$, then either
$(j,k)$, $(k,l)$ and $(l,s)$ are three different pairs or
$(j,k)=(l,s)\neq (k,l)$. In any case Lemma~\ref{alp4a} implies that
$M_{j,k,l,s}=\{0\}$ for generic $d$-dimensional $L_{j,k}$ ($j\neq
k$). According to the last display $\E_4(L,E)=\{0\}$ for generic
$d$-dimensional $L_{j,k}$ ($j\neq k$). Thus, there exists a
$d(m^2-m)$-dimensional subspace $L$ of $E\otimes E$ such that
$\E_4(L,E)=\{0\}$. By Lemma~\ref{tepr},
$h_4(\K,m^2,m^4-(m^2-m)d)=0$. Since $n\leq m^2$, we get
$h_4(\K,n,m^4-(m^2-m)d)=0$. Hence
$$
d(\K,n,4)\leq m^4-(m^2-m)d.
$$
Since $d=m^2-r^2$ and $r<\frac{\sqrt5-1}{2}m+1$, we have
$d>\frac{\sqrt 5-1}2 m^2-(\sqrt 5-1)m-1$. Using this inequality
together with $\sqrt n\leq m\leq \sqrt n+1$ and the fact that the
functions $m\mapsto m^2-m$ and $m\mapsto \frac{\sqrt 5-1}2
m^2-(\sqrt 5-1)m-1$ on $[1,\infty)$ are increasing, we see that the
above display implies
\begin{align*}
d(\K,n,4)&\leq (\sqrt n +1)^4-(n-\sqrt n)\Bigl(\frac{\sqrt 5-1}2
n-(\sqrt 5-1)\sqrt n -1\Bigr)=
\\
&=\frac{3-\sqrt 5}{2}n^2+\frac{5+3\sqrt 5}{2}n^{3/2}+(8-\sqrt 5)
n+3\sqrt n+1.
\end{align*}
The above display immediately implies (\ref{ine}).
\end{proof}

\subsection{An estimate of $d(\K,n,4)$ for any field}

We prove the following specific lemma in the appendix.

\begin{lemma}\label{3-4} Let $\K$ be any field and
$$
R=\K\langle x_1,x_2,x_3\rangle/{\tt
Id}\{x_1x_2,x_1x_3,x_2x_3,x_1^2+x_2^2+x_3^2\}.
$$
Then the Hilbert series of $R$ is $1+3t+5t^2+4t^3$. In particular,
$h_4(\K,3,4)=0$ for any field $\K$.
\end{lemma}

\begin{corollary}\label{3-41}
Let $\K$ be a field, $E$ be a three-dimensional vector space over
$\K$ with a basis $\{e_1,e_2,e_3\}$ and $L$ be the $5$-dimensional
subspace of $E\otimes E$ spanned by $e_2\ootimes e_1$, $e_3\ootimes
e_2$, $e_3\ootimes e_1$, $e_1\ootimes e_1-e_2\ootimes e_2$ and
$e_1\ootimes e_1-e_3\ootimes e_3$. Then $\E_4(L,E)=\{0\}$, where
$\E_k(L,E)$ are defined in $(\ref{ekl})$.
\end{corollary}

\begin{proof} Let $R$ be the algebra defined in Lemma~\ref{3-4}.
From Lemmas~\ref{3-4} and~\ref{rel} it follows that $0=\dim
R_4=\dim\E_4(M^\perp,E)$, where $M$ is the 4-dimensional subspace of
$E\otimes E$ spanned by the vectors $e_1\ootimes e_2$, $e_1\ootimes
e_3$, $e_2\ootimes e_3$ and $e_1\ootimes e_1+e_2\ootimes
e_2+e_3\ootimes e_3$. It is straightforward to verify that
$L=M^\perp$. Thus $\E_4(L,E)=\{0\}$.
\end{proof}

\begin{proposition}\label{gfield} Let $\K$ be a field. Then the
numbers $d(\K,n,4)$ defined in $(\ref{dnq})$ satisfy the following
inequality
\begin{equation}\label{gg}
d(\K,n,4)\leq d_n=\left\{\begin{array}{ll}\frac{4n^2}9&\text{if $n=3k,$ $k\in\N;$}\\
\frac{4n^2+2n-2}9&\text{if $n=3k+2,$ $k\in\Z_+;$}\\
\frac{4n^2+4n-8}9&\text{if $n=3k+1,$ $k\in\N$.}
\end{array} \right.
\end{equation}
In any case $d(\K,n,4)<\frac{4(n^2+n)}9$.
\end{proposition}

\begin{proof} Let $E_3$ be a 3-dimensional vector space over $\K$ with
a basis $\{e_1,e_2,e_3\}$ and $E_2=\spann\{e_1,e_2\}$. Consider the
$5$-dimensional subspace $L$ of $E_3\otimes E_3$ spanned by
$e_2\ootimes e_1$, $e_3\ootimes e_2$, $e_3\ootimes e_1$,
$e_1\ootimes e_1-e_2\ootimes e_2$ and $e_1\ootimes e_1-e_3\ootimes
e_3$ and let $L_{2,2}=(E_2\otimes E_2)\cap L$, $L_{2,3}=(E_2\otimes
E_3)\cap L$ and $L_{3,2}=(E_3\otimes E_2)\cap L$. By
Corollary~\ref{3-41}, $\E_4(L,E_3)=\{0\}$. It is also easy to see
that $L_{2,3}=L_{2,2}$ is spanned by $e_2\ootimes e_1$ and
$e_1\ootimes e_1-e_2\ootimes e_2$ and therefore $\dim L_{2,3}=\dim
L_{2,2}=2$. On the other hand $L_{3,2}$ is spanned by $e_2\ootimes
e_1$, $e_3\ootimes e_2$, $e_3\ootimes e_1$ and $e_1\ootimes
e_1-e_2\ootimes e_2$ and therefore $\dim L_{3,2}=4$.

Let $n=3k$ with $k\in\N$. Then the direct sum $G$ of $k$ copies of
$E_3$ has dimension $n$. By Lemma~\ref{main} with
$n_1={\dots}=n_k=3$, there exists a subspace $L_G$ of $G\otimes G$
of dimension $k^2\dim L=5k^2=\frac{5n^2}9$ such that
$\E_4(L_G,G)=\{0\}$. By Lemma~\ref{tepr},
$h_4(\K,n,\frac{4n^2}9)=0$. Hence $d(\K,n,4)\leq \frac{4n^2}9$.

Let $n=3k+2$ with $k\in\Z_+$. Then the direct sum $G$ of $E_2$ and
$k$ copies of $E_3$ has dimension $n$. By Lemma~\ref{main} with
$n_1=2$ and $n_2={\dots}=n_{k+1}=3$, there is a subspace $L_G$ of
$G\otimes G$ such that $\E_4(L_G,G)=\{0\}$ and $\dim L_G=k^2\dim
L+k(\dim L_{2,3}+\dim L_{3,2})+\dim
L_{2,2}=5k^2+6k+2=\frac{5n^2-2n+2}9$. By Lemma~\ref{tepr},
$h_4(\K,n,\frac{4n^2+2n-2}9)=0$. Hence $d(\K,n,4)\leq
\frac{4n^2+2n-2}9$.

Let $n=3k+1$ with $k\in\N$. Then the direct sum $G$ of $2$ copies of
$E_2$ and $k-1$ copies of $E_3$ has dimension $n$. By
Lemma~\ref{main} with $n_1=n_2=2$ and $n_3={\dots}=n_{k+1}=3$, there
is a subspace $L_G$ of $G\otimes G$ such that $\E_4(L_G,G)=\{0\}$
and $\dim L_G=(k-1)^2\dim L+2(k-1)(\dim L_{2,3}+\dim L_{3,2})+4\dim
L_{2,2}=5k^2+2k+1= \frac{5n^2-4n+8}9$. By Lemma~\ref{tepr},
$h_4(\K,n,\frac{4n^2+4n-8}9)=0$. Hence $d(\K,n,4)\leq
\frac{4n^2+4n-8}9$.
\end{proof}

\begin{corollary}\label{co2} Let $\K$ be any field. Then
$h_4(\K,n,\frac{n(n-1)}2)=0$ for any $n\geq 17$ and for
$n\in\{9,12,14,15\}$.
\end{corollary}

\begin{proof} Let $A=\{9,12,14,15\}\cup\{n\in\N:n\geq 17\}$.
By Proposition~\ref{gfield}, $h_4(\K,n,d_n)=0$, where $d_n$ is
defined in (\ref{gg}). It is straightforward to verify that
$\frac{n(n-1)}{2}\geq d_n$ for $n\in A$. Hence
$h_4(\K,n,\frac{n(n-1)}2)=0$ for $n\in A$.
\end{proof}

\subsection{Proof of Theorem \ref{main0}}

Proposition \ref{gfield} was the main step in the proof of Theorem \ref{main0}. It ensures that for $d \geq \frac{4(n^2+n)}{9}$ the fourth component of the generic Hilbert series vanishes. Now we combine this with the already known due to Anick fact (Theorem A), that the third component of the generic series always coincides with the third component of the Golod-Shafarevich series (which vanishes for $d \geq \frac{n^2}{2}$). So, by looking at the third and fourth components of the Hilbert series we obtain the statement of Theorem \ref{main0}.

\subsection{$\alpha(\K,5)$ for any field}

The proof of the following Lemma is quite technical and is presented
in the appendix.

\begin{lemma}\label{3-3}Let $\K$ be a field and $R=\K\langle
x_1,x_2,x_3\rangle/I$ with $I={\tt Id}\{f_1,f_2,f_3\}$ and
\begin{equation*}
f_1=x_3^2-x_1x_2,\ \ f_2=x_3x_2-x_2x_3+x_2x_1-x_1x_3-x_1x_2+x_1^2,\
\ f_3=x_3x_1+x_2^2-x_1^2.
\end{equation*}
Then the Hilbert series of $R$ is
$H_R(t)=1+3t+6t^2+9t^3+9t^4=|(1-3t+3t^2)^{-1}|$.
\end{lemma}

\begin{corollary}\label{3-31}
Let $\K$ be a field, $E$ be a three-dimensional vector space over
$\K$ with a basis $\{e_1,e_2,e_3\}$ and $L$ be the $5$-dimensional
subspace of $E\otimes E$ spanned by Consider the $6$-dimensional
subspace $L$ of $E_3\otimes E_3$ spanned by $e_3\ootimes
e_3+e_1\ootimes e_2+e_1\ootimes e_1+e_2\ootimes e_2$, $e_2\ootimes
e_3+e_1\ootimes e_1+e_2\ootimes e_2$, $e_1\ootimes e_3+e_1\ootimes
e_1+e_2\ootimes e_2$, $e_2\ootimes e_1-e_1\ootimes e_1-e_2\ootimes
e_2$, $e_3\ootimes e_2-e_1\ootimes e_1-e_2\ootimes e_2$ and
$e_3\ootimes e_1+e_1\ootimes e_1-e_2\ootimes e_1$. Then
$\E_5(L,E)=\{0\}$.
\end{corollary}

\begin{proof} Let $R$ be the algebra defined in Lemma~\ref{3-3}.
From Lemmas~\ref{3-3} and~\ref{rel} it follows that $0=\dim
R_5=\dim\E_5(M^\perp,E)$, where $M$ is the 3-dimensional subspace of
$E\otimes E$ spanned by the vectors $e_3\ootimes e_3-e_1\ootimes
e_2$, $e_3\ootimes e_2-e_2\ootimes e_3+e_2\ootimes e_1-e_1\ootimes
e_3-e_1\ootimes e_2+e_1\ootimes e_1$ and $e_3\ootimes
e_1+e_2\ootimes e_2-e_1\ootimes e_1$. It is straightforward to
verify that $L=M^\perp$. Thus $\E_5(L,E)=\{0\}$.
\end{proof}

\begin{proposition}\label{whatnot1} Let $\K$ be any field. Then the
numbers $d(\K,n,5)$ defined in $(\ref{dnq})$ satisfy the following
inequality
\begin{equation}\label{gg5}
d(\K,n,5)\leq \delta_n=\left\{\begin{array}{ll}\frac{n^2}3&\text{if $n=3k,$ $k\in\N;$}\\
\frac{n^2+2n+1}3&\text{if $n=3k+2,$ $k\in\Z_+;$}\\
\frac{n^2+3n+1}3&\text{if $n=3k+1,$ $k\in\N$.}
\end{array} \right.
\end{equation}
In particular, $\alpha(\K,5)=\frac13$ for any field $\K$.
\end{proposition}

\begin{proof}Let $E_3$ be a 3-dimensional vector space over $\K$ with
a basis $\{e_1,e_2,e_3\}$, $E_2=\spann\{e_1,e_2\}$ and
$E_1=\spann\{e_1\}$. Consider the $6$-dimensional subspace $L$ of
$E_3\otimes E_3$ defined in Corollary~\ref{3-31} and let
$L_{j,k}=(E_j\otimes E_k)\cap L$ for $1\leq j,k\leq 3$. By
Corollary~\ref{3-31}, $\E_5(L,E_3)=\{0\}$. Estimating from below the
dimensions of $L_{j,k}$ by the number of basic vectors of $L$
contained in $E_j\otimes E_k$, we get $\dim L_{3,3}=\dim L=6$, $\dim
L_{2,3}\geq 3$, $\dim L_{3,2}\geq 3$, $\dim L_{2,2}\geq 1$ and $\dim
L_{3,1}\geq 1$.

Let $n=3k$ with $k\in\N$. Then the direct sum $G$ of $k$ copies of
$E_3$ has dimension $n$. By Lemma~\ref{main} with
$n_1={\dots}=n_k=3$, there exists a subspace $L_G$ of $G\otimes G$
of dimension $k^2\dim L=6k^2=\frac{2n^2}3$ such that
$\E_5(L_G,G)=\{0\}$. By Lemma~\ref{tepr}, $h_5(\K,n,\frac{n^2}3)=0$.
Hence $d(\K,n,4)\leq \frac{n^2}3$.

Let $n=3k+2$ with $k\in\Z_+$. Then the direct sum $G$ of $E_2$ and
of $k$ copies of $E_3$ has dimension $n$. By Lemma~\ref{main} with
$n_1=2$ and $n_2={\dots}=n_{k+1}=3$, there exists a subspace $L_G$
of $G\otimes G$ such that $\E_5(L_G,G)=\{0\}$ and $\dim L_G=k^2\dim
L+k(\dim L_{2,3}+\dim L_{3,2})+\dim L_{2,2}\geq
6k^2+6k+1=\frac{2n^2-2n-1}3$. By Lemma~\ref{tepr},
$h_4(\K,n,\frac{n^2+2n+1}3)=0$. Hence $d(\K,n,4)\leq
\frac{n^2+2n+1}3$.

Let $n=3k+1$ with $k\in\N$. Then the direct sum $G$ of $E_1$ and $k$
copies of $E_3$ has dimension $n$. By Lemma~\ref{main} with $n_1=1$
and $n_2={\dots}=n_{k+1}=3$, there exists a subspace $L_G$ of
$G\otimes G$ such that $\E_5(L_G,G)=\{0\}$ and $\dim L_G=k^2\dim
L+k(\dim L_{1,3}+\dim L_{3,1})+\dim L_{1,1}\geq 6k^2+k=
\frac{2n^2-3n-1}3$. By Lemma~\ref{tepr},
$h_5(\K,n,\frac{n^2+3n+1}3)=0$. Hence $d(\K,n,4)\leq
\frac{n^2+3n+1}3$.

Thus we have verified (\ref{gg5}), from which it follows that
$\alpha(\K,5)\leq\frac13$. On the other hand, for $0\leq
d<\frac{n^2}{3}$, the coefficient in $|(1-nt+dt^2)^{-1}|$ in front
of $t^5$ is $n^5-4n^3d+3nd^2>0$. By Theorem~GS, $h_5(\K,n,d)>0$ for
$d<\frac{n^2}3$. It follows that $\alpha(\K,5)\geq \frac13$. Thus
$\alpha(\K,5)=\frac13$.
\end{proof}

\section{Further applications of Lemma~\ref{main}}

As we have seen in the last section, specific examples with 3
generators produce non-trivial estimates for $d(\K,n,4)$ and
$d(\K,n,5)$ for any $\K$ and $n$. We proceed along the same lines.
To this end we need more specific examples of Hilbert series of
quadratic algebras. We produce them via lucky guesswork and
application of the software package GRAAL ('Graded Algebras' by A.Verevkin and
A.Kondratiev) to find the Hilbert series. This program uses the
one-sided Gr\"obner basis technique to calculate the Hilbert series
of $\Z_p$-algebras given by generators and relations.

\begin{example}\label{7-19} Let $R=\Z_2\langle
x_1,\dots,x_7\rangle/I$ with $I={\tt Id}\{f_1,\dots,f_{19}\}$ and
\begin{equation*}
\begin{array}{lll}
f_1=x_1x_7,& f_2=x_3x_7+x_4x_6+x_6x_2,&
f_3=x_5x_7+x_6x_4+x_3x_5+x_2x_1+x_4x_3,\\
f_4=x_7x_1+x_1x_6,& f_5=x_7x_2+x_6x_1+x_1x_5,&
f_6=x_1^2+x_2^2+x_3^2+x_4^2+x_5^2+x_6^2+x_7^2,\\
f_7=x_2x_7+x_7x_3,& f_8=x_6x_7+x_3x_6+x_4x_5+x_5x_2,&
f_9=x_7x_5+x_2x_6+x_5x_3+x_1x_4+x_3x_2,\\
f_{10}=x_5x_7+x_7x_6,& f_{11}=x_7x_6+x_6x_2+x_5x_1+x_3x_4,&
f_{12}=x_7x_4+x_6x_3+x_2x_5+x_3x_2+x_4x_1,\\
f_{13}=x_7x_5+x_4x_7,& f_{14}=x_7x_2+x_3x_6+x_6x_4,&
f_{15}=x_2x_7+x_6x_5+x_5x_4+x_3x_1+x_4x_2,\\
f_{16}=x_3x_7+x_7x_4,& f_{17}=x_4x_7+x_2x_6+x_6x_3,&
f_{18}=x_7x_3+x_4x_6+x_5x_2+x_2x_4+x_3x_1,\\
\rlap{\hbox{$f_{19}=x_6x_7+x_6x_4+x_2x_6+x_2x_5+x_3x_5+x_4x_5$.}}
\end{array}
\end{equation*}
Then the Hilbert series of $R$ is
$H_R(t)=1+7t+30t^2+77t^3=|(1-7t+19t^2)^{-1}|$.
\end{example}

\begin{example}\label{4-6}Let $R=\Z_2\langle
x_1,x_2,x_3,x_4\rangle/I$ with $I={\tt Id}\{f_1,\dots,f_{6}\}$ and
\begin{equation*}
\begin{array}{lll}
f_1=x_1x_2,& f_2=x_1x_4+x_4x_2+x_2x_3,&
f_3=x_1x_3+x_3x_4+x_4x_1,\\
f_4=x_1^2+x_2^2+x_3^2+x_4^2,& f_5=x_3x_4+x_4x_2+x_2x_4,&
f_6=x_2x_3+x_3x_1+x_1x_3.
\end{array}
\end{equation*}
Then the Hilbert series of $R$ is
$H_R(t)=1+4t+10t^2+16t^3+t^4=|(1-4t+6t^2)^{-1}|$.
\end{example}

\begin{example}\label{4-5}Let $R=\Z_2\langle
x_1,x_2,x_3,x_4\rangle/I$ with $I={\tt Id}\{f_1,\dots,f_{5}\}$ and
\begin{equation*}
\begin{array}{lll}
f_1=x_1^2+x_2^2+x_3^2+x_4^2,& f_2=x_1x_2+x_2x_3+x_3x_4,&
f_3=x_4x_1+x_1x_3+x_3x_2,\\
f_4=x_1x_3+x_3x_2+x_2x_4,& f_5=x_1x_4+x_4x_3+x_3x_2+x_2x_4.&
\end{array}
\end{equation*}
Then the Hilbert series of $R$ is
$H_R(t)=1+4t+11t^2+24t^3+41t^4+44t^5=|(1-4t+5t^2)^{-1}|$.
\end{example}

\subsection{Proof of Theorem~\ref{lith1}}

The equalities $\alpha(\K,4)=\frac{3-\sqrt5}{2}$ and
$\alpha(\K,5)=\frac13$ follow from Propositions~\ref{whatnot}
and~\ref{whatnot1} respectively.

Let $\K$ be a field of characteristic $0$. Example~\ref{4-5} shows
that $d(\Z_2,4,6)\leq 5$. By Proposition~\ref{fiel} $d(\K,4,6)\geq
5$. By Lemma~\ref{limi}
$$
\alpha(\K,\infty)\leq\alpha(\K,6)=\lim_{n\to\infty}\frac{d(\K,n,6)}{n^2}=
\inf_{n\in\N}\frac{d(\K,n,6)}{n^2}\leq\frac{d(\K,4,6)}{16}\leq\frac{5}{16}.
$$
This completes the proof of Theorem~\ref{lith1}.

\subsection{Proof of Theorem~\ref{solv}}

Throughout this section $\K$ is a field of characteristic 0, $n\geq
3$. Lemma~\ref{3-3} together with Theorem~GS show that $H^{\rm
min}_{\K,3,3}(t)=1+3t+6t^2+9t^3+9t^4=|(1+3t-3t^2)^{-1}|$, which
proves Theorem~\ref{solv} for $n=3$.

Example~\ref{4-6} and Proposition~\ref{fiel} show that $H^{\rm
min}_{\K,4,6}(t)\leq H^{\rm min}_{\Z_2,4,6}(t)\leq
1+4t+10t^2+16t^3+t^4=|(1-4t+6t^2)^{-1}|$. By Theorem~GS, $H^{\rm
min}_{\K,4,6}(t)\geq |(1-4t+6t^2)^{-1}|$. Hence $H^{\rm
min}_{\K,4,6}(t)= 1+4t+10t^2+16t^3+t^4=|(1-4t+6t^2)^{-1}|$, which
proves Theorem~\ref{solv} for $n=4$.

From now on $n\geq 5$. It suffices to prove that
\begin{equation}\label{44}
h_4(\K,n,\delta_n)=0,\ \ \ \text{where $\delta_n=\frac{n(n-1)}{2}$.}
\end{equation}
Indeed, Theorem~A implies that $h_3(\K,n,\delta_n)=n^2$ and
therefore $H^{\rm
min}_{\K,n,\delta_n}(t)=1+nt+\frac{n(n+1)}2t^2+n^2t^3=|(1+nt-\delta_nt^2)^{-1}|$
provided (\ref{44}) is true. Thus it remains to prove (\ref{44}) for
$n\geq 5$.

If $n\in A=\{9,12,14,15\}\cup\{n\in\N:n\geq 17\}$,
Corollary~\ref{co2} implies that (\ref{44}) is satisfied. It remains
to consider $n\in\{5,6,7,8,10,11,13,16\}$.

According to the remarks in the beginning of Section~\ref{genob},
the relations in Example~\ref{7-19} considered as relations in
$\K\langle x_1,\dots,x_7\rangle$ define a quadratic $\K$-algebra $R$
such that $H_R(t)\leq 1+7t+30t^2+77t^3=|(1-7t+19t^2)^{-1}|$. On the
other hand, by Theorem~GS, $H_R(t)\geq |(1-7t+19t^2)^{-1}|$. Hence
$$
H_R(t)=1+7t+30t^2+77t^3.
$$
Let $E$ be a $7$-dimensional vector space over $\K$ with a basis
$\{e_1,\dots,e_7\}$ and $E_j=\spann\{e_1,\dots,e_j\}$ for $1\leq
j\leq 7$. Let also $M$ be the $19$-dimensional subspace of $E\otimes
E$ spanned by the elements obtained from the relations of $R$ by
replacing $x_jx_k$ by $e_j\ootimes e_k$. It is straightforward to
verify that the $30$-dimensional subspace $L=M^\perp$ of $E\otimes
E$ is spanned by the following $30$ vectors: $g_1=e_1\ootimes e_2$,
$g_2=e_2\ootimes e_3$, $g_3=e_1\ootimes e_3$, $g_4=e_1\ootimes
e_1-e_2\ootimes e_2$, $g_5=e_1\ootimes e_1-e_3\ootimes e_3$,
$g_6=e_2\ootimes e_4-e_3\ootimes e_1+e_4\ootimes e_2$,
$g_7=e_1\ootimes e_4-e_3\ootimes e_2+e_4\ootimes e_1$,
$g_8=e_1\ootimes e_1-e_4\ootimes e_4$, $g_9=e_2\ootimes
e_1-e_4\ootimes e_3$, $g_{10}=e_1\ootimes e_1-e_5\ootimes e_5$,
$g_{11}=e_5\ootimes e_1-e_3\ootimes e_4$, $g_{12}=e_5\ootimes
e_4-e_4\ootimes e_2$, $g_{13}=e_5\ootimes e_3-e_3\ootimes
e_2+e_4\ootimes e_1$, $g_{14}=e_5\ootimes e_2-e_4\ootimes
e_5-e_2\ootimes e_4+e_3\ootimes e_5$, $g_{15}=e_2\ootimes
e_5-e_3\ootimes e_5-e_3\ootimes e_2+e_1\ootimes e_4$,
$g_{16}=e_1\ootimes e_1-e_6\ootimes e_6$, $g_{17}=e_6\ootimes
e_1-e_1\ootimes e_5$, $g_{18}=e_5\ootimes e_6$, $g_{19}=e_4\ootimes
e_6-e_6\ootimes e_2+e_3\ootimes e_4-e_2\ootimes e_4$,
$g_{20}=e_6\ootimes e_4-e_2\ootimes e_1-e_3\ootimes e_5-e_3\ootimes
e_6+e_5\ootimes e_2-e_2\ootimes e_4$, $g_{21}=e_2\ootimes
e_6-e_1\ootimes e_4-e_6\ootimes e_3-e_3\ootimes e_5+e_4\ootimes
e_1$, $g_{22}=e_6\ootimes e_5-e_4\ootimes e_2$, $g_{23}=e_1\ootimes
e_1-e_7\ootimes e_7$, $g_{24}=e_7\ootimes e_1-e_1\ootimes e_6$,
$g_{25}=e_3\ootimes e_7-e_7\ootimes e_4-e_4\ootimes e_6+e_2\ootimes
e_4+e_4\ootimes e_1$, $g_{26}=e_5\ootimes e_7-e_7\ootimes
e_6-e_2\ootimes e_1+e_3\ootimes e_4$, $g_{27}=e_7\ootimes
e_2-e_6\ootimes e_1-e_6\ootimes e_4+e_3\ootimes e_5$,
$g_{28}=e_2\ootimes e_7-e_7\ootimes e_3+e_2\ootimes e_4-e_4\ootimes
e_2$, $g_{29}=e_6\ootimes e_7-e_3\ootimes e_5-e_5\ootimes
e_2+e_2\ootimes e_4$ and $g_{30}=e_7\ootimes e_5-e_4\ootimes
e_7+e_6\ootimes e_3-e_4\ootimes e_1-e_1\ootimes e_4$.

As usual, $L_{j,k}=L\cap(E_j\otimes E_k)$ for $2\leq j,k\leq 7$.
Denoting $d_j=\dim L_{j,j}$ for $1\leq j\leq 7$, we easily obtain
$d_7=30$, $d_6=22$, $d_5=15$ and $d_4=9$. Denoting $d_{j,k}=\dim
L_{j,k}+\dim L_{k,j}$ for $1\leq j<k\leq 7$ and estimating $\dim
L_{l,m}$ from below by the number of basic vectors $g_s$ in
$E_l\otimes E_m$, we get $\dim d_{6,7}\geq 47$, $d_{5,6}\geq 33$ and
$d_{4,6}\geq 22$.

By Lemma~\ref{rel}, the equality $R_4=\{0\}$ implies
\begin{equation}\label{e4}
\E_4(L,E)=0.
\end{equation}
Condition (\ref{rel}) implies that $\E_4(L_{n,n},E_n)=\{0\}$ for
$1\leq n\leq 7$. Using this equality and applying Lemma~\ref{tepr},
we get $h_4(\K,n,n^2-d_n)=0$ for $1\leq n\leq 7$. This equality for
$n\in\{5,6,7\}$ gives $h_4(\K,5,10)=h_4(\K,6,14)=h_4(\K,7,19)=0$.
Since $\delta_5=10$, $\delta_6=15>14$ and $\delta_7=21>19$,
(\ref{44}) is satisfied for $n\in\{5,6,7\}$.

Using Lemma~\ref{main} with $m=2$ and $n_1=n_2=4$, we see that there
exists a subspace $N$ of $G\otimes G$ such that $\dim G=8$, $\dim
N\geq 4d_4=36$ and $\E_4(N,G)=0$. By Lemma~\ref{tepr},
$h_4(\K,8,28)=0$, which is (\ref{44}) for $n=8$.

Using Lemma~\ref{main} with $m=2$ and $n_1=n_2=5$, we see that there
exists a subspace $N$ of $G\otimes G$ such that $\dim G=10$, $\dim
N\geq 4d_5=60$ and $\E_4(N,G)=0$. By Lemma~\ref{tepr},
$h_4(\K,10,40)=0$. Hence $h_4(\K,10,45)=0$ which is (\ref{44}) for
$n=10$.

Using Lemma~\ref{main} with $m=2$ and $n_1=5$ and $n_2=6$, we see
that there exists a subspace $N$ of $G\otimes G$ such that $\dim
G=11$, $\dim N\geq d_6+d_5+d_{5,6}\geq 22+15+33=70$ and
$\E_4(N,G)=0$. By Lemma~\ref{tepr}, $h_4(\K,11,51)=0$. Hence
$h_4(\K,11,55)=0$ which is (\ref{44}) for $n=11$.

Using Lemma~\ref{main} with $m=2$ and $n_1=6$ and $n_2=7$, we see
that there exists a subspace $N$ of $G\otimes G$ such that $\dim
G=13$, $\dim N\geq d_7+d_6+d_{6,7}\geq 30+22+47=99$ and
$\E_4(N,G)=0$. By Lemma~\ref{tepr}, $h_4(\K,13,70)=0$. Hence
$h_4(\K,13,78)=0$ which is (\ref{44}) for $n=13$.

Using Lemma~\ref{main} with $m=3$ and $n_1=n_2=6$ and $n_3=4$, we
see that there exists a subspace $N$ of $G\otimes G$ such that $\dim
G=16$, $\dim N\geq 4d_6+d_4+2d_{4,6}\geq 4\cdot 22+9+2\cdot 22=141$
and $\E_4(N,G)=0$. By Lemma~\ref{tepr}, $h_4(\K,16,115)=0$. Hence
$h_4(\K,16,120)=0$ which is (\ref{44}) for $n=16$.

The proof of Theorem~\ref{solv} is complete.

\small

\section{Appendix: Proof of Lemmas~\ref{3-4} and~\ref{3-3}}

In the proofs of both Lemmas~\ref{3-4} and~\ref{3-3} we use the
non-commutative analog of the Buchberger algorithm of constructing a
Gr\"obner basis.

\subsection{Proof of Lemma~\ref{3-4}}

Ordering the variables as $x_1>x_2>x_3$ and considering the
degree-lexicographical ordering on the monomials, one can easily see
that the set
$$
\{x_2x_3,x_1,x_3,x_1x_2,x_1^2+x_2^2+x_3^2,
x_2^3+x_3^2x_2,x_2^2x_1+x_3^2x_1,x_3^2,x_2^2,x_3^2x_2x_1\}
$$
is the reduced Gr\"obner basis of the ideal $I={\tt
Id}\{x_1x_2,x_1x_3,x_2x_3,x_1^2+x_2^2+x_3^2\}$. By analyzing the
above basis, one can easily verify that
$\{x_3^2x_2+I,x_3^2x_1+I,x_3x_2^2+I,x_3x_2x_1+I\}$ is a linear basis
in $R_3$ and that $R_4=\{0\}$. Hence the Hilbert series of $R$ is
$1+3t+5t^2+4t^3$. The proof is complete.

\subsection{Proof of Lemma~\ref{3-3}}

We consider the ordering $x_1>x_2>x_3$ on the variables and the
corresponding degree-lexicographical ordering on the monomials.

First, we assume that the characteristic of $\K$ is different from
$2$. Using the non-commutative analog of the Buchberger algorithm of
constructing a Gr\"obner basis, we find the following homogeneous
elements of $I$, written starting from the leading term (=highest
monomial) and having the property that neither of the leading terms
are subwords of the others:

\begin{align*}
g_1&=x_1x_3{-}x_2x_1{-}x_2^2{+}x_2x_3{-}x_3x_1{-}x_3x_2{+}x_3^2,\ \
g_2=x_1x_2{-}x_3^2,\ \
g_3=x_1^2{-}x_1x_3{+}x_2x_1{-}x_2x_3{+}x_3x_2{-}x_3^2,
\\
g_4&=x_2^2x_3{+}x_2x_3x_1{+}2x_2x_3x_2{-}3x_2x_3^2{-}x_3x_1^2{+}2x_3x_1x_3{-}x_3x_2^2{+}x_3x_2x_3{-}x_3^2x_1{+}x_3^2x_2{-}x_3^3,
\\
g_5&=x_2^3{+}x_3x_2x_1{-}x_3^2x_1{+}x_3^2x_2,\quad
g_6=x_2^2x_1{+}x_2x_3x_1{+}x_2x_3x_2{-}2x_2x_3^2{-}x_3x_1x_2{+}x_3x_1x_3{+}x_3x_2x_3{-}x_3^3,
\\
g_7&=2x_2x_3x_2x_1{+}x_3x_1^3{-}x_3x_1^2x_3{+}6x_3x_1x_2x_1{-}3x_3x_1x_2^2{-}x_3x_1x_2x_3{-}x_3x_1x_3x_1{+}3x_3x_1x_3x_2{-}
\\
&\quad
{-}2x_3x_1x_3^2{-}5x_3x_2x_1^2{-}4x_3x_2x_1x_2{+}6x_3x_2x_1x_3{-}2x_3x_2^2x_1{+}x_3x_2^3{+}2x_3x_2^2x_3{-}3x_3x_2x_3x_1{+}
\\
&\quad
{+}x_3x_2x_3x_2{+}3x_3x_2x_3^2{+}3x_32x_1^2{+}4x_3^2x_1x_2{+}x_3^2x_1x_3{-}8x_3^2x_2x_1{-}9x_3^2x_2^2{+}9x_3^2x_2x_3{-}5x_3^3x_2{-}3x_3^4,
\\
g_8&=x_2x_3x_2^2{-}x_3x_1^3{-}x_3x_1^2x_3{-}3x_3x_1x_2x_1{+}3x_3x_1x_2^2{-}x_3x_1x_2x_3{+}2x_3x_1x_3x_1{+}x_3x_1x_3x_2{-}x_3x_1x_3^2{+}2x_3x_2x_1^2{-}
\\
&\quad{-}3x_3x_2x_1x_2{-}4x_3x_2x_1x_3{+}2x_3x_2^2x_1{-}x_3x_2^2x_3{+}x_3x_2x_3x_1{+}x_3x_2x_3x_2{+}x_3x_2x_3^2{-}
\\
&\quad
{-}3x_3^2x_1^2{+}3x_3^2x_1x_2{+}2x_3^2x_1x_3{+}x_3^2x_2x_1{+}4x_3^2x_2^2{-}4x_3^2x_2x_3{+}3x_3^3x_1{-}x_3^3x_2{+}2x_3^4,
\\
g_9&=2x_2x_3x_2x_3{+}x_3x_1^3{+}9x_3x_1x_2x_1{-}5x_3x_1x_2^2{+}2x_3x_1x_2x_3{-}x_3x_1x_3x_1{+}2x_3x_1x_3x_2{-}2x_3x_1x_32{-}10x_3x_2x_1^2{-}
\\
&\quad{-}5x_3x_2x_1x_2{+}18x_3x_2x_1x_3{-}4x_3x_2^2x_1{-}2x_3x_2^3{+}4x_3x_2^2x_3{-}4x_3x_2x_3x_1{-}4x_3x_2x_3x_2{+}3x_3x_2x_3^2{+}
\\
&\quad{+}7x_3^2x_1^2{+}5x_3^2x_1x_2{-}4x_3^2x_1x_3{-}16x_3^2x_2x_1{-}19x_3^2x_2^2{+}24x_3^2x_2x_3{-}4x_3^3x_1{+}4x_3^3x_2{-}13x_3^4,
\\
g_{10}&=2x_2x_3^2x_1{+}x_3x_1^3{-}x_3x_1^2x_3{+}4x_3x_1x_2x_1{-}x_3x_1x_2^2{-}x_3x_1x_2x_3{-}x_3x_1x_3x_1{+}3x_3x_1x_3x_2{-}2x_3x_1x_3^2{-}5x_3x_2x_1^2{-}
\\
&\quad{-}6x_3x_2x_1x_2{+}6x_3x_2x_1x_3{+}x_3x_2^3{+}2x_3x_2^2x_3{-}3x_3x_2x_3x_1{+}x_3x_2x_3x_2{+}3x_3x_2x_3^2{+}3x_3^2x_1^2{+}
\\
&\quad{+}6x_3^2x_1x_2{+}x_3^2x_1x_3{-}10x_3^2x_2x_1{-}9x_3^2x_2^2{+}9x_3^2x_2x_3{+}2x_3^3x_1{-}5x_3^3x_2{-}3x_3^4,
\\
g_{11}&=2x_2x_3^2x_2{-}x_3x_1^3{-}3x_3x_1x_2x_1{+}3x_3x_1x_2^2{+}x_3x_1x_3x_1{-}5x_3x_2x_1x_2{+}2x_3x_2^2x_1{+}x_3x_2x_32{-}x_3^2x_1^2{+}
\\
&\quad{+}5x_3^2x_1x_2{-}2x_3^2x_2x_1{-}x_3^2x_2^2{+}4x_3^3x_1{-}x_3^4,
\\
g_{12}&=x_2x_3^3{+}x_3x_1x_2x_1{+}x_3x_1x_3x_2{-}x_3x_1x_3^2{-}2x_3x_2x_1^2{-}3x_3x_2x_1x_2{+}3x_3x_2x_1x_3{+}x_3x_2^2x_3{-}x_3x_2x_3x_1{+}x_3x_2x_3^2{+}
\\
&\quad{+}x_3^2x_1^2{+}3x_3^2x_1x_2{-}4x_3^2x_2x_1{-}4x_3^2x_2^2{+}4x_3^2x_2x_3{+}x_3^3x_1{-}x_3^3x_2{-}2x_3^4,
\\
g_{13}&=x_3^2x_2x_3x_1,\quad g_{14}=x_3^2x_2x_3x_2,\quad
g_{15}=x_3^2x_2x_3^2,\quad g_{16}=x_3^3x_2x_1,\quad
g_{17}=x_3^3x_2^2,
\\
g_{18}&=x_3^3x_2x_3,\quad g_{19}=x_3^4x_1,\quad g_{20}=x_3^4x_2,
\quad g_{21}=x_3^5.
\end{align*}

It follows that the dimension of $R_q$ does not exceed the number of
monomials of degree $q$ that do not contain a subword being the
leading monomial of one of $g_j$. This observation gives $\dim
R_3\leq 9$, $\dim R_4\leq 9$ and $\dim R_5=0$. Thus $H_R(t)\leq
1+3t+6t^2+9t^3+9t^4=|(1-3t+3t^2)^{-1}|$. On the other hand,
Theorem~GS implies that $H_R(t)\geq |(1-3t+3t^2)^{-1}|$. Hence
$H_R(t)=1+3t+6t^2+9t^3+9t^4=|(1-3t+3t^2)^{-1}|$.

It remains to consider the case ${\rm char}\,\K=2$. Using the same
algorithm, we find the following homogeneous elements of $I$,
written starting from the leading term (=highest monomial) and
having the property that neither of the leading terms are subwords
of the others:
\begin{align*}
g_1&=x_1^2+x_2^2+x_3x_1,\quad g_2=x_1x_2+x_3^2,\quad g_3=
x_1x_3+x_2x_1+x_2^2+x_2x_3+x_3x_1+x_3x_2+x_3^2,
\\
g_{4}&=x_2^2x_1+x_2x_3x_1+x_2x_3x_2+x_3x_2x_1+x_3x_2^2+x_3^2x_1+x_3^2x_2+x_3^3,
\\
g_{5}&=x_2^3+x_3x_2x_1+x_3^2x_1+x_3^2x_2,\quad
g_{6}=x_2^2x_3+x_2x_3x_1+x_2x_3^2+x_3x_2x_3+x_3^2x_2+x_3^3,
\\
g_7&=x_2x_3x_2x_1+x_2x_3^2x_1+x_3x_2x_3x_2+x_3^2x_2x_1+x_3^3x_1,
\\
g_8&=x_2x_3x_2^2+x_2x_3^2x_2+x_2x_3^3+x_3x_2x_3x_2+x_3x_2x_3^2+x_3^2x_2x_1+x_3^2x_2x_3+x_3^3x_2+x_3^4,
\\
g_9&=x_2x_3x_2x_3+x_2x_3^3+x_3x_2x_3x_1+x_3^2x_2^2+x_3^2x_2x_3+x_3^3x_1+x_3^4,
\\
g_{10}&=x_2x_3^2x_2+x_3x_2x_3x_2+x_3x_2x_3^2+x_3^2x_2x_1+x_3^2x_2^2+x_3^3x_1+x_3^4,
\\
g_{11}&=x_2x_3^3+x_3x_2x_3x_2+x_3x_2x_3^2+x_3^2x_2x_1+x_3^3x_1,\quad
g_{12}=x_3x_2x_3x_1+x_3x_2x_3^2+x_3^2x_2x_1+x_3^2x_2^2,
\\
g_{13}&=x_3x_2x_3^2x_1,\quad g_{14}=x_3^2x_2x_3x_2,\quad
g_{15}=x_3^2x_2x_3^2,\quad g_{16}=x_3^3x_2x_1,
\\
g_{17}&= x_3^3x_2^2, \quad g_{18}= x_3^3x_2x_3,\quad
g_{19}=x_3^4x_1,\quad g_{20}= x_3^4x_2,\quad g_{21}= x_3^5.
\end{align*}

Again, the dimension of $R_q$ does not exceed the number of
monomials of degree $q$ that do not contain a subword being the
leading monomial of one of $g_j$. This observation gives $\dim
R_3\leq 9$, $\dim R_4\leq 9$ and $\dim R_5=0$. Thus $H_R(t)\leq
1+3t+6t^2+9t^3+9t^4=|(1-3t+3t^2)^{-1}|$. On the other hand,
Theorem~GS implies that $H_R(t)\geq |(1-3t+3t^2)^{-1}|$. Hence
$H_R(t)=1+3t+6t^2+9t^3+9t^4=|(1-3t+3t^2)^{-1}|$. The proof is
complete.

\normalsize


\rm

\small

\vskip1truecm

\scshape



\noindent Queens's University Belfast

\noindent Department of Pure Mathematics

\noindent University road, Belfast, BT7 1NN, UK

\noindent E-mail: \qquad {\tt n.iyudu@qub.ac.uk}\ \ \ {\rm
and}\ \ \ {\tt s.shkarin@qub.ac.uk}

\end{document}